\documentclass{article}
\usepackage{amsmath,amssymb,amsthm, amscd, multicol, mathrsfs, tikz, graphics, float, hyperref}

\newtheorem{definition}[subsection]{Definition}
\newtheorem{theorem}[subsection]{Theorem}
\newtheorem{corollary}[subsection]{Corollary}

\newtheorem{proposition}[subsection]{Proposition}

\newtheorem{observation}[subsection]{Observation}
\renewcommand{\qed}{\hfill$\Box$\\}

\begin{document}
\title{Minimum rank and failed zero forcing number of graphs}
\author{Prince Allan B. Pelayo\thanks{%
Institute of Mathematics, University of the Philippines, Diliman, Quezon
City 1101; \texttt{papelayo@math.upd.edu.ph}} \; and Ma. Nerissa M. Abara \thanks{%
Institute of Mathematics, University of the Philippines, Diliman, Quezon
City 1101; \texttt{issa@math.upd.edu.ph}}}
\maketitle

\begin{abstract}

Let $G$ be a simple, finite, and undirected graph with vertices each given an initial coloring of either blue or white. Zero forcing on graph $G$ is an iterative process of forcing its white vertices to become blue after a finite application of a specified color-change rule. We say that an initial set $S$ of blue vertices of $G$ is a zero forcing set for $G$ under the specified color-change rule if a finite number of iterations of zero forcing results to an updated coloring where all vertices of $G$ are blue. Otherwise, we say that $S$ is a failed zero forcing set for $G$ under the specified color-change rule. It is not difficult to see that any subset of a failed zero forcing set is also failed. Hence, our interest lies on the maximum possible cardinality of a failed zero forcing set, which we refer to as the failed zero forcing number of $G$. In this paper, we consider two color-change rules $-$ standard and positive semidefinite. We compute for the failed zero forcing numbers of several graph families. Furthermore, under each graph family, we characterize the graphs $G$ for which the failed zero forcing number is equal to the minimum rank of $G$. \newline
\bigskip


\medskip

 \end{abstract}

\vspace{0.3in}

\section{Introduction}\label{sec1}

A \textit{graph} is an ordered pair $G = (V, E)$ comprised of a non-empty set $V$ of vertices and a (possibly empty) set $E$ of edges. We say that $G$ is a \textit{simple, finite, and undirected graph} if it satisfies all of the following: (a) each edge is a pair of two distinct vertices (that is, $G$ has no loops); (b) each pair of vertices comprises at most one edge (that is, $G$ has no multiple edges); (c) $V$ is finite and; (d) the edges are unordered pairs of vertices. \\
\indent Throughout this paper, we consider simple, finite, and undirected graphs. For terminology and notation not defined here, we refer to \cite{W}. We denote the \textit{vertex set} and \textit{edge set} of a graph $G$ by $V(G)$ and $E(G)$, respectively. The cardinality of $V(G)$ is called the \textit{order of $G$}, denoted $|G|$. A vertex $v$ adjacent to a vertex $u$ is said to be a \textit{neighbor of $u$}. A vertex without any neighbor is called an \textit{isolated vertex}. A set $K \subset V(G)$ of $k$ vertices having the same neighbors outside $K$ is called a \textit{module of order $k$}. Two non-adjacent vertices in a module are said to be \textit{similar}. Finally, we use the following notations in graph theory:
\begin{center}
\begin{tabular}{lll}
$N(u)$ & \phantom{white} & set of neighbors of ve
rtex $u$ \\
deg$(u)$ & \phantom{white} & degree of vertex $u$ \\
$\Delta(G)$ & \phantom{white} & maximum degree of graph $G$ \\
$\delta(G)$ & \phantom{white} & minimum degree of graph $G$ \\
$G[S]$ & \phantom{white} & subgraph of $G$ induced by set $S$ \\
$P_n$  & \phantom{white} & path with $n$ vertices \\
$C_n$  & \phantom{white} & cycle with $n$ vertices \\
$K_n$  & \phantom{white} & complete graph with $n$ vertices \\
$W_n$  & \phantom{white} & wheel with $n$ vertices \\
$K_{m,n}$  & \phantom{white} & complete bipartite graph on $m$ and $n$ vertices \\
$Q_n$ & \phantom{white} & $n$th hypercube \\
$H_s$ & \phantom{white} & $s$th half-graph \\
$\overline{K_n}$ & \phantom{white} & graph with $n$ isolated vertices \\
\end{tabular}
\end{center}

\indent There is a correspondence between graphs and matrices. In particular, given a graph $G$, there is an infinite number of matrices whose graph is $G$. Let $F$ be a field. Consider the set of $n$-by-$n$ symmetric matrices whose entires are in $F$, denoted ${S}_n(F)$. The graph of a matrix $A$ in ${S}_n(F)$, denoted $\mathscr{G}(A)$, is the graph $G$ whose vertex set and edge set are given by $V(G) = \{1,  \ldots  , n\}$ and $E(G) = \{ij : a_{ij} \neq 0 \; \text{and} \; 1 \leq i < j \leq n\}$, respectively. Note that the diagonal entries of $A$ are ignored when determining the graph of $A$. \\
\indent The set of symmetric matrices of graph $G$ over $F$ is defined to be the set:
$$\mathscr{S}(F,G) = \{A \in S_n(F) : \mathscr{G}(A) = G\}$$
\indent On the other hand, the set of positive semidefinite matrices of graph $G$ over $F$ is defined to be the set: 
$$\mathscr{S}_+(F,G) = \{A \in \mathscr{S}(F,G) : A \; \text{is positive semidefinite}\}$$
\indent In this paper, we consider $F = \mathbb{R}$. For simplicity, we use the notations $\mathscr{S}(G)$ and $\mathscr{S}_+(G)$ to denote $\mathscr{S}(\mathbb{R}, G)$ and $\mathscr{S}_+(\mathbb{R}, G)$, respectively. The minimum rank and maximum nullity of a graph $G$ (over $\mathbb{R}$) is defined, respectively, to be:
$$\text{mr}(G) = \text{min}\{\text{rank}(A) : A \in \mathscr{S}(G)\}$$
$$\text{M}(G) = \text{max}\{\text{null}(A) : A \in \mathscr{S}(G)\}$$
\indent On the other hand, the minimum positive semidefinite rank and maximum positive semidefinite nullity of a graph $G$ (over $\mathbb{R}$) is defined, respectively, to be:
$$\text{mr}_+(G) = \text{min}\{\text{rank}(A) : A \in \mathscr{S}_+(G)\}$$
$$\text{M}_+(G) = \text{max}\{\text{null}(A) : A \in \mathscr{S}_+(G)\}$$
\indent The \textit{minimum rank problem} for a graph $G$ is a problem which aims to find the value of $\text{mr}(G)$. There is also an interest in a variant of this problem $-$ the \textit{minimum positive semidefinite rank problem}. In this variant, the value of $\text{mr}_+(G)$ is computed instead of $\text{mr}(G)$. By the rank-nullity theorem, the minimum (positive semidefinite) rank problem is also equivalent to determining the maximum (positive semidefinite) nullity of a graph. We refer to \cite{FH} for a survey of known results on minimum rank problems. \\
\indent The significant amount of attention given to the minimum rank problem has paved way for the introduction of several related graph parameters. One such parameter is the zero forcing number, which was introduced in \cite{AIM}. A chapter in \cite{FH2} is dedicated to a thorough discussion on the relationship between zero forcing and minimum rank of graphs. On the other hand, the failed-type parameter of the zero forcing number was considered in \cite{FJS}. The notion of the failed zero forcing number of a graph is a relatively new topic and thus poses a lot of open problems. \\
\indent This paper is organized as follows. In Section \ref{sec2}, we introduce failed zero forcing sets for a graph and present some properties of such sets. We also introduce the graph parameter $F(G)$ and $F_+(G)$, the maximum cardinality of a failed (positive semidefinite) zero forcing set. In Sections \ref{sec3} and \ref{sec4}, we compute values of the failed (positive semidefinite) zero forcing number for several graph families. Finally, in Section \ref{sec5}, we relate $F(G)$ and $F_+(G)$ with other graph parameters including minimum rank and maximum nullity.

\section{Failed zero forcing sets and failed positive \\ semidefinite zero forcing sets}\label{sec2}

We begin with a graph $G$ with each vertex of $G$ having an initial coloring of either blue or white. We let $S$ to be the set of all blue vertices. Given a color-change rule, the \textit{derived coloring} is the result of iteratively applying the color-change rule on $G$ with $S$ as the initial set until no more changes are possible. 

\begin{definition}
\textit{Standard zero forcing} 
\begin{itemize}
\item \textit{Standard color-change rule:} {\normalfont If $v \in V(G) \backslash S$ is the unique white vertex of $G$ that is adjacent to a blue vertex $u$ in $S$, then change the color of $v$ from white to blue.}
\item \textit{Iteration:} {\normalfont Update $S$ by adjoining all vertices made blue.}
\item {\normalfont We say that a subset $S$ of $V(G)$ is a \textit{zero forcing set} for $G$ if the derived coloring under the standard color-change rule with $S$ as the initial set is all blue. Otherwise, we say that $S$ is a \textit{failed zero forcing set.}}
\item {\normalfont The \textit{zero forcing number of a graph $G$}, denoted $Z(G)$, is the minimum of $|Z|$ taken over all zero forcing sets $Z$ for $G$.}
\item {\normalfont The \textit{failed zero forcing number of a graph $G$}, denoted $F(G)$, is the maximum of $|F|$ taken over all failed zero forcing sets $F$ for $G$.}
\end{itemize}
\end{definition}

\begin{definition}
\textit{Positive semidefinite zero forcing} 
\begin{itemize}
\item \textit{Positive semidefinite color-change rule:} {\normalfont Identify the sets $W_1,  \ldots , W_k$ of the $k$ connected components of $G[V(G) \backslash S]$. Consider each $i \in \{1, \ldots, k\}$. If $v \in W_i$ is the unique white vertex of $G[W_i \cup S]$ that is adjacent to a blue vertex in $S$, then change the color of $v$ from white to blue.}
\item \textit{Iteration:} {\normalfont Update $S$ by adjoining all vertices made blue.}
\item {\normalfont We say that a subset $S$ of $V(G)$ is a \textit{positive semidefinite zero forcing set} for $G$ if the derived coloring under the positive semidefinite color-change rule with $S$ as the initial set is all blue. Otherwise, we say that $S$ is a \textit{failed positive semidefinite zero forcing set.}}
\item {\normalfont The \textit{positive semidefinite zero forcing number of a graph $G$}, denoted $Z_+(G)$, is the minimum of $|Z|$ taken over all positive semidefinite zero forcing sets $Z$ for $G$.}
\item {\normalfont The \textit{failed positive semidefinite zero forcing number of a graph $G$}, denoted $F_+(G)$, is the maximum of $|F|$ taken over all failed positive semidefinite zero forcing sets $F$ for $G$.}
\end{itemize}
\end{definition}
\indent Note that the positive semidefinite color-change rule can be restated in the following manner: If $W_1,  \ldots , W_k$ are the $k$ connected components arising from the white vertices of $G$, then apply the standard color-change rule to $S$ in $G[W_i \cup S]$ for all $i \in \{1,  \ldots , k\}$ and update $S$ by adjoining all vertices made blue. \\
\indent If a blue vertex $u$ changes the color of a vertex $v$ from white to blue, we say that $u$ forces $v$, denoted $u \rightarrow v$. 

\begin{observation}\label{obs2.3}
Let $S$ be a zero forcing set (resp. positive semidefinite zero forcing set) for a graph $G$. If $S'$ is a subset of $V(G)$ that contains $S$, then $S'$ is also a zero forcing set (resp. positive semidefinite zero forcing set) for $G$. 
\end{observation}

\begin{observation}\label{obs2.4}
Let $S$ be a failed zero forcing set (resp. failed positive semidefinite zero forcing set) for a graph $G$. If $S'$ is a subset of $V(G)$ that is contained in $S$, then $S'$ is also a failed zero forcing set (resp. failed positive semidefinite zero forcing set) for $G$. 
\end{observation}

\indent Observations \ref{obs2.3} and \ref{obs2.4} justify why it is natural and intuitive to consider the minimum size of a (positive semidefinite) zero forcing set and the maximum size of a failed (positive semidefinite) zero forcing set. A (positive semidefinite) zero forcing set with cardinality equal to the (positive semidefinite) zero forcing number is said to be a \textit{minimum (positive semidefinite) zero forcing set}. On the other hand, a failed (positive semidefinite) zero forcing set with cardinality equal to the failed (positive semidefinite) zero forcing number is said to be a \textit{maximum failed (positive semidefinite) zero forcing set}. \\
\indent We also consider \textit{minimal (positive semidefinite) zero forcing sets} and \textit{maximal failed (positive semidefinite) zero forcing sets}. Given a property $\mathscr{P}$, a set $P$ is said to be minimal with respect to $\mathscr{P}$ if $P$ satisfies $\mathscr{P}$ but any $P'$ properly contained in $P$ does not satisfy $\mathscr{P}$. On other hand, the set $P$ is said to be maximal with respect to $\mathscr{P}$ if $P$ satisfies $\mathscr{P}$ but any $P'$ properly containing $P$ does not satisfy $\mathscr{P}$.

\begin{definition}\label{def2.5} 
Let $G$ be a graph. Given a color-change rule $\mathscr{R}$, we say that a proper subset $S$ of $V(G)$ is \textit{stalled} under $\mathscr{R}$ if no color changes are possible from $S$. 
\end{definition}

\begin{observation}\label{obs2.6}
Maximal and maximum failed (positive semidefinite) zero forcing sets are stalled under the (positive semidefinite) color-change rule. Furthermore, a maximum failed (positive semidefinite) zero forcing set is maximal. 
\end{observation}

\indent In \cite{AIM}, it is shown that the elements of a zero forcing set can be interpreted as a set of indices associated with a coordinate in a vector that is required to be zero, while a white vertex corresponds to a coordinate that can be either zero or nonzero. Given a graph $G$ and a vector in the kernel of a matrix in $\mathscr{S}(G)$, the term 'zero forcing' pertains to the idea that all coordinates must be zero (that is, 'forced' to be zero) whenever all coordinates associated to the blue vertices are initially set to zero. Let $x = [x_i]$ be a vector in $F^{n \times 1}$, where $F$ is a field and $n>0$. The \textit{support} of $x$, denoted $supp(x)$, is the set of indices $i$ such that $x_i \neq 0$. The next two propositions are lifted from \cite{AIM}. From these, we give a necessary and sufficient condition for a set to be a failed zero forcing set. 

\begin{proposition}\label{prop2.7} 
\cite[Proposition $2.2$]{AIM} Suppose $F$ is a field, $A \in F^{n \times n}$, and  $null(A) > k$. If $W$ is a set of $k$ indices, then there is a nonzero vector $x \in ker(A)$ such that $supp(x) \cap W = \emptyset$. 
\end{proposition}

\begin{proposition}\label{prop2.8} 
\cite[Proposition $2.3$]{AIM}  Let $Z$ be a zero forcing set for a  graph $G$ and $A \in \mathscr{S}(G)$. If $x \in ker(A)$ and $supp(x) \cap Z = \emptyset$, then $x=\vec{0}$. 
\end{proposition}

Since any subset of $V(G)$ with $F(G) + 1$ elements is a zero forcing set for a graph $G$, the next corollary is an immediate consequence of Proposition \ref{prop2.8}.

\begin{corollary}\label{cor2.9}
Let $S \subseteq V(G)$ with $|S|$ $\geq$ $F(G) + 1$, and let $A$ $\in$ $\mathscr{S}(G)$. If $x \in ker(A)$ and all entries of $x$ corresponding to $S$ are zero, then $x = \vec{0}$. 
\end{corollary}

\begin{corollary}\label{cor2.10}
Let $S \subseteq V(G)$ and $A$ $\in$ $\mathscr{S}(G)$. If there exists a nonzero vector $x \in ker(A)$ such that $supp(x)$ $\cap$ $S$ $=$ $\emptyset$, then $S$ is a failed zero forcing set for $G$. If $null(A) > F(G)$, then the converse also holds. 
\end{corollary}
\proof
Let $x \in \text{ker}(A) \backslash \{\vec{0}\}$ such that $\text{supp}(x) \cap S = \emptyset$. If $S$ is a zero forcing set for $G$, then $x = \vec{0}$ by Proposition \ref{prop2.8}, which is a contradiction. Hence, $S$ is failed. Now, suppose $null(A) > F(G)$. Let $S$ be a failed zero forcing set. Then, $F(G) \geq |S|$. This implies that $null(A) \geq |S|$. By Proposition \ref{prop2.7}, there exists a nonzero vector $x \in ker(A)$ such that $supp(x) \cap S$ $=$ $\emptyset$. 
\qed

Analogous results for positive semidefinite zero forcing set (as consequences of \cite[Theorem $3.5$]{Barioli}) are given by the next two propositions. 

\begin{proposition}\label{prop2.11}
Let $S \subseteq V(G)$ with $|S|$ $\geq$ $F_+(G) + 1$, and let $A$ $\in$ $\mathscr{S}_+(G)$. If  $x \in ker(A)$ and all entries of $x$ corresponding to $S$ are zero, then $x = \vec{0}$. 
\end{proposition}

\begin{proposition}\label{prop2.12}
Let $S \subseteq V(G)$ and $A$ $\in$ $\mathscr{S}_+(G)$. If there exists a nonzero vector $x \in ker(A)$ such that $supp(x)$ $\cap$ $S$ $=$ $\emptyset$, then $S$ is a failed positive semidefinite zero forcing set for $G$. If $null(A) > F_+(G)$, then the converse also holds. 
\end{proposition}

\section{The failed zero forcing number}\label{sec3}

\begin{table}[t]
\begin{center}
\begin{tabular}{llllll}
\hline
Result & $G$ & $F(G)$ & $F(G) =$ mr$(G)$? \\
\hline		
\ref{thm3.6}, \ref{thm5.7} & $P_n$ & $\lceil \frac{n-2}{2} \rceil$ & iff $n=1$ \\ 
\ref{thm3.6}, \ref{thm5.7} & $C_n$, $n \geq 3$ & $\lfloor \frac{n}{2} \rfloor$ & iff $n=3, 4$ \\
\ref{thm3.6}, \ref{thm5.7} &$K_n$, $n \geq 2$ & $n-2$ & iff $n=3$ \\
\ref{thm3.6}, \ref{thm5.7} & $W_4$ & 2 & no  \\
\ref{thm3.6}, \ref{thm5.7} & $W_5$ & 3 & no \\
\ref{thm3.6}, \ref{thm5.7} & $W_n, n \geq 6$ & $\lfloor \frac{2n-2}{3} \rfloor$ & iff $n=6, 7$ \\
\ref{thm3.6}, \ref{thm5.7} & $K_{m,1}$, $m \geq 1$ & $m-1$ & iff $m =3$ \\
\ref{thm3.6}, \ref{thm5.7} & $K_{m,2}$, $m \geq 2$ & $m$ & iff $m=2$ \\
\ref{thm3.6}, \ref{thm5.7} & $K_{m,n}, m \geq n \geq 2$ & $m+n-2$ & iff $m+n=4$ \\
\ref{thm3.7}, \ref{thm5.7} & $Q_1$ & $0$ & no  \\
\ref{thm3.7}, \ref{thm5.7} & $Q_2$ & $2$ & yes  \\
\ref{thm3.7}, \ref{thm5.7} & $Q_n, n \geq 3$ & $\geq 2^{n-1} - n$ & no \\
\ref{thm3.8}, \ref{thm5.7} &$H_1$ & 0 & no \\
\ref{thm3.8}, \ref{thm5.7} & $H_s, s \geq 2$ & $2s-3$ & iff $s=3$ \\
\hline
\end{tabular}
\caption{\label{table6.1} Summary of results on $F(G)$ in this paper}
\end{center}
\end{table}

The graph parameter $F(G)$ was introduced in $\cite{FJS}$. In this section, we enumerate some results from $\cite{FJS}$ while extending the study to consider other graph families. 

\begin{observation}\label{obs3.1}
\cite[Observation $2.1$]{FJS} For any graph $G$, $Z(G) - 1 \leq F(G) \leq |G| - 1$. 
\end{observation}

\begin{proposition}\label{prop3.2}
Let $G$ be a disconnected graph with $G_1,  \ldots , G_k$ as its connected components. Suppose $i \in \{1, \ldots , k\}$. If
$$F = V(G) \big\backslash \big( V(G_i) \backslash F_i \big),$$
where $F_i$ is a maximum failed zero forcing set for $G_i$, then $F$ is a maximal failed zero forcing set for $G$. 
\end{proposition}
\proof
Suppose $G$ is disconnected with connected components $G_1,  \ldots , G_k$. Let $i \in \{1,  \ldots , k\}$. Let $F = V(G) \backslash (V(G_i) \backslash F_i)$, where $F_i$ is a maximum failed zero forcing set for $G_i$. Since $F_i$ is failed for $G_i$, then $V(G_i) \backslash F_i) \neq \emptyset$. This further implies that $F \neq V(G)$. Let $v \in V(G) \backslash F$. Consider $u \in F$. Then $u \in F_i$ or $u \in V(G_j)$ for some $j \neq i$. If $u \in F_i$, then $u$ cannot force $v$ since $F_i$ is failed for $G_i$. If $u \in V(G_j)$ for some $j \neq i$, then $u$ is non-adjacent to $v$ since $G_i$ and $G_j$ are disconnected from each other. Hence, no vertex can force $v$. Thus, $F$ is stalled. Since $F \neq V(G)$, then $F$ is a maximal failed zero forcing set for $G$. 
\qed

\begin{corollary}\label{cor3.3}
Let $G$ be a disconnected graph with $G_1,  \ldots , G_k$ as its connected components. Then $F(G) = |G| -$ min$\{|G_i| - F(G_i)$ $:$ $i = 1,  \ldots , k\}$.
\end{corollary}

\proof
Suppose $G$ is a disconnected graph. Let $G_1,  \ldots , G_k$ be its connected components. Suppose $F_i$ is a maximum failed zero forcing set for $G_i$ for all $i \in \{1,  \ldots , k\}$. By Proposition \ref{prop3.2}, $F(G) \geq |V(G) \backslash (V(G_i) \backslash F_i)|$ $=$ $|G| - (|G_i|-F(G_i))$. Hence, $F(G) \geq |G| -$ min$\{|G_i| - F(G_i)$ $:$ $i = 1,  \ldots , k\}$. Suppose $S$ is a subset of $V(G)$ such that $|S| > |G| -$ min$\{|G_i| - F(G_i)$ $:$ $i = 1,  \ldots , k\}$. Then, $|S| > |G| - (|G_i| - F(G_i))$ for $i = 1,  \ldots , k$. Manipulating the inequality, we get $|S| + |G_i| - |G| > F(G_i)$. Note that $|S \cap V(G_i)|$ $=$ $|S| + |G_i| - |S \cup V(G_i)|$ $\geq$ $|S| + |G_i| - |G|$ $>$ $F(G_i)$. Since $F(G_i)$ is the maximum value of any failed set for $G_i$, we must have that the set $S \cap V(G_i)$ is a zero forcing set for $G_i$. Since this is true for all $i \in \{1,  \ldots , k\}$, $S$ is also a zero forcing set for $G$. Hence, $F(G) \leq |G| -$ min$\{|G_i| - F(G_i)$ $:$ $i = 1,  \ldots , k\}$. Therefore, $F(G) = |G| -$ min$\{|G_i| - F(G_i)$ $:$ $i = 1,  \ldots , k\}$. 
\qed

\begin{observation}\label{obs3.4}
\cite[Observation $2.1$]{FJS} Suppose $G$ is a graph with $n$ vertices. Then $F(G) = n -1$ if and only if $G$ has an isolated vertex. 
\end{observation}

\begin{theorem}\label{thm3.5}
\cite[Theorem $2.4$]{FJS} Suppose $G$ is a connected graph with $n$ vertices. Then $F(G) = n-2$ if and only if $G$ has a module of order 2. 
\end{theorem}

\begin{theorem}\label{thm3.6}
\cite[Theorems $3.1$, $3.2$, $3.3$, $3.4$, $3.6$, and $3.8$]{FJS} Let $G$ be a simple, finite, and connected graph. Then
$$F(G) =
\begin{cases}
\lceil \frac{n-2}{2} \rceil, & \text{if} \hspace{2mm} G = P_n \\ 
\lfloor \frac{n}{2} \rfloor, & \text{if} \hspace{2mm} G = C_n, n \geq 3 \\
n - 2, & \text{if} \hspace{2mm} G = K_n, n \geq 2 \\
n - 2, & \text{if} \hspace{2mm} G \hspace{1.7mm} \text{is an m-ary tree on n vertices.} \hspace{1.7mm} m \geq 2 \\
3, & \text{if} \hspace{2mm} G = W_n, n = 5\\
\lfloor \frac{2n-2}{3} \rfloor, & \text{if} \hspace{2mm} G = W_n, n \geq 4, n \neq 5 \\
m+n-2 & \text{if} \hspace{2mm} G = K_{m,n}, m+n \geq 3, m \geq n 
\end{cases}$$
\end{theorem}

\begin{theorem}\label{thm3.7}
Let $G$ be the hypercube $Q_n$. Then $F(Q_1) = 0$, $F(Q_2) = 2$, and $F(Q_n) \geq 2^n - n$ for $n \geq 3$.
\end{theorem}
\proof
Note that $Q_1 = P_2$ and $Q_2 = C_4$. By Theorem \ref{thm3.6}, $F(Q_1) = 0$ and $F(Q_2) = 2$. Suppose $n \geq 3$ and let $G = Q_n$. Label the $2^n$ vertices of $G$ using binary $n$-tuples. That is, if $v \in V(G)$, then we can write $v$ as $(x_1, x_2, \ldots , x_n)$, where $x_i$ is either 0 or 1 for each $i \in \{1, \ldots, n\}$. In this labeling, two vertices $(x_1, x_2, \ldots, x_n)$ and $(y_1, y_2, \ldots, y_n)$ are adjacent if and only if $x_i = y_i$ except for exactly one value of $i \in \{1, \ldots, n\}$. \\
\indent Consider the set $F = V(G) \backslash W$, where $W$ is the set of all vertices for which $x_i = 1$ except for exactly one value of $i \in \{1,  \ldots , n\}$. Let a vertex $v$ be blue if $v \in F$. Otherwise, let $v$ be white. \\
\indent Let $v$ be a blue vertex in $F$. Then, $v$ has one of the following forms: (1) $v = (1, 1, \ldots, 1)$; (2) $v = (x_1, x_2, \ldots, x_n)$ with at least three $x_i$'s equal to 0; (3) $v = (x_1, x_2, \ldots, x_n)$ with exactly two $x_i$'s equal to 0. \\
\indent If $v = (1, 1, \ldots , 1)$, then $v$ is adjacent to every vertex in $W$. Since $n \geq 3$, then $|W| \geq 3$. Hence, $v$ is adjacent to at least three white vertices. \\
\indent If $v = (x_1, x_2, \ldots , x_n)$ where $m$ of the $x_i$'s are equal to 0 with $m \geq 3$, then $v$ can only be adjacent to those vertices with $m-1$ (which is at least 2) or $m+1$ (which is at least 4) of the $x_i$'s equal to 0. Such vertices are already blue. \\
\indent Suppose $v = (x_1, x_2, \ldots , x_n)$ has exactly two $x_i$'s equal to 0, say $x_j = x_k = 0$ with $1 \leq j < k \leq n$. Then the vertex $v$ is adjacent to $w = (y_1, y_2, \ldots, y_n)$ , where $y_l = 0$ if and only if $l = j$, and to $u = (z_1, z_2, \ldots, z_n)$, where $z_l=0$ if and only if $l=k$. Hence $w, u \in W$ and are white. \\
\indent Thus, every vertex in $F$ with a neighbor outside $F$ is adjacent to at least two white vertices. Thus, $F$ is stalled and failed. Finally, we have that $F(Q_n) \geq |F| = |Q_n| - |W| = 2^n - {{n}\choose{n-1}} = 2^n - n$ for $n \geq 3$.
\qed

\begin{theorem}\label{thm3.8}
Let $G$ be the half-graph $H_s$. Then $F(H_1) = 0$ and $F(H_s) = 2s - 3$ for $s \geq 2$. 
\end{theorem}
\proof
Note that $H_1 = P_2$ and $H_2 = P_4$. By Theorem \ref{thm3.6}, $F(H_1) = 0$ and $F(H_2) = 1$. Suppose $s > 2$. Since the half-graph is a conencted graph with no isolated vertices nor a module of order $2$, then $F(H_s) \leq 2s - 3$. Next, partition $H_s$ into $S_1$ and $S_2$ such that this partition defines the bipartite property of $H_s$. Label the vertices in $S_1$ as $v_1$,  \ldots , $v_s$ and the corresponding vertices in $S_2$ as $v_{s+1},  \ldots , v_{2s}$, respectively. Let $V = \{v_{s}, v_{2s-1}, v_{2s}\}$. Note that $v_{s}$ cannot be forced since its only neighbor $v_{2s}$ is outside $V(H_s)$ $\setminus$ $V$. On the other hand, every vertex in $V(H_s)$ $\setminus$ $V$ adjacent to $v_{2s-1}$ is also adjacent to $v_{2s}$. Therefore, $v_{2s-1}$ and $v_{2s}$ cannot be forced. This implies that the set $V(H_s)$ $\setminus$ $V$ is stalled. Therefore, $V(H_s)$ $\setminus$ $V$ is a failed zero forcing set and $F(H_s) \geq 2s - 3$. Hence, $F(H_s) = 2s - 3$. 
\qed

\section{The failed positive semidefinite zero forcing number}\label{sec4}

\begin{table}[t]
\begin{center}
\begin{tabular}{llllll}
\hline
Result & $G$ & $F_+(G)$ & $F_+(G) =$ mr$_+(G)$? \\
\hline		
\ref{thm4.5}, \ref{thm5.8} & $P_n$ & 0 & iff $n=1$ \\ 
\ref{thm4.6}, \ref{thm5.8} & $C_n$, $n \geq 3$ & 1 & iff $n=3$ \\
\ref{cor4.13}, \ref{thm5.8} & $K_n$, $n \geq 2$ & $n-2$ & $n=3$ \\
\ref{thm4.20}, \ref{thm5.8} & $W_4$ & 2 & no  \\
\ref{thm4.20}, \ref{thm5.8} & $W_5$ & 2 & yes \\
\ref{thm4.20}, \ref{thm5.8} & $W_n, n \geq 6$ & $\lfloor \frac{2n-2}{3} \rfloor$ & iff $n=5, 6,7$ \\
\ref{thm4.21}, \ref{thm5.8} & $K_{m,1}$, $m \geq 1$ & 0 & no \\
\ref{thm4.21}, \ref{thm5.8} & $K_{m,2}$, $m \geq 2$ & $m-1$ & no \\
\ref{thm4.21}, \ref{thm5.8} & $K_{m,n}, m \geq n \geq 2$ & $m+n-4$ & iff $n=4$ \\
\ref{thm4.22}, \ref{thm5.8} & $Q_1$ & $0$ & no \\
\ref{thm4.22}, \ref{thm5.8} & $Q_2$ & 1 & no \\
\ref{thm4.22}, \ref{thm5.8} & $Q_n, n \geq 3$ & $\geq 2^{n-1} - n - 1$ & iff $n=3$ \\
\ref{thm4.23}, \ref{thm5.8} & $H_1$ & 0 & no \\
\ref{thm4.23}, \ref{thm5.8} & $H_s, s \geq 2$ & $2s-4$ & iff $s=4$ \\
\hline
\end{tabular}
\caption{\label{table6.2} Summary of results on $F_+(G)$ in this paper}
\end{center}
\end{table}

Both failed-type parameter of standard zero forcing \cite{FJS} and skew zero forcing \cite{Ansill} were already previously introduced. In this paper, we formally introduce and study the failed-type parameter of positive semidefinite zero forcing. 

\indent We begin by characterizing graphs for which $F_+(G)$ is very close to $|G|$. The vertex set $V(G)$ of a graph $G$ is trivially a positive semidefinite zero forcing set since there are no other vertices left to be forced. Also, any set with cardinality less than $Z_+(G)$ is not a positive semidefinite zero forcing set. Hence, we have the following proposition:

\begin{proposition}\label{prop4.1}
For any graph $G$, $Z_+(G) - 1 \leq F_{+}(G) \leq |G|-1$. 
\end{proposition}

\begin{theorem}\label{thm4.2}
If $G$ is a graph with $n$ vertices, then $F_{+}(G) = n-1$ if and only if $G$ has an isolated vertex. 
\end{theorem}
\proof
Let $G$ be a graph with an isolated vertex $v$. From \ref{prop4.1}, $F_+(G) \geq n-1$. Equality holds since $V(G)$ $\backslash$ $\{v\}$ is failed. Suppose $F_{+}(G) = n-1$ and suppose $G$ has no isolated vertex. Consider a failed positive semidefinite zero forcing set $S$ with $|S| = n-1$. Let $v$ be the vertex of $G$ not in $S$. Since $G$ has no isolated vertex, there exists $u \in S$ such that $uv \in E(G)$. Hence, $u \rightarrow v$ and $S$ is a positive semidefinite zero forcing set for $G$, which is a contradiction to our assumption. Hence, $G$ has an isolated vertex. 
\qed

\begin{proposition}\label{prop4.3}
Let $G$ be a disconnected graph with $n$ vertices. If one of the connected components of $G$ is a path with $k$ vertices, where $k > 1$, then $F_+(G) \geq n-k$. If one of the connected components of $G$ is a cycle with $m$ vertices, where $m > 2$, then $F_+(G) \geq n - (m-1)$. 
\end{proposition}
\proof Let $G$ be a disconnected graph. Suppose a path $P_k$ is one of the connected components of $G$, where $k > 1$. Color all vertices in the set $S = V(G) \backslash V(P_k)$ blue. Consider a vertex $v \notin S$. Then, $v \in V(P_k)$. Since all vertices in $V(P_k)$ are not in $S$ and $P_k$ is disconnected from the rest of the graph, no vertex in $S$ can force $v$. Therefore, $S$ is stalled and is a failed positive semidefinite zero forcing set. Hence, $F_+(G) \geq |S| = |V(G) \backslash V(P_k)| = n - k$. On the other hand, if one of the connected components of $G$ is a cycle $C_m$, where $m > 2$, color all vertices in the set $T$ blue, where $T = ( V(G) \backslash V(C_m) ) \cup \{v\}$ for some $v \in V(C_m)$. Since $m > 2$, $v$ has two neighbors outside $T$. The rest of the vertices in $T$ cannot force a vertex outside $T$ since $C_m$ is disconnected from the rest of the graph. Thus, $T$ is stalled and is failed. Hence, we have $F_+(G) \geq |T| =$ $|(V(G) \backslash V(C_m)) \cup \{v\}| = (n - m) + 1 = n - (m-1)$.
\qed

\begin{theorem}\label{thm4.5}
If $G$ is a path, then $F_+(G) = 0$. 
\end{theorem}
\proof
The empty set is trivially a failed positive semidefinite zero forcing set for the graph. Hence, $F_+(P_n) \geq 0$. Label the vertices starting from an endpoint as $v_1, \ldots, v_n$. Consider the set $\{v_i\}$ for some $i \in \{1, \ldots, n\}$ and color $v_i$ blue. If $i=1$, then $v_1 \rightarrow v_2 \rightarrow  \ldots  \rightarrow v_n$. If $i = n$, then $v_n \rightarrow v_{n-1} \rightarrow  \ldots  \rightarrow v_1$. Finally, if $i \in \{2,  \ldots , n-1\}$, then $\{v_1,  \ldots , v_{i-1}\}$ and $\{v_{i+1},  \ldots , v_n\}$ are two connected components arising from the white vertices having $v_i$ as an endpoint. Hence, $v_i \rightarrow v_{i-1} \rightarrow  \ldots  \rightarrow v_1$ and $v_i \rightarrow v_{i+1} \rightarrow  \ldots  \rightarrow v_n$. This implies that $S$ is a positive semidefinite zero forcing set and $F_+(P_n) < 1$. Therefore, $F_+(P_n) = 0$. 
\qed

\begin{theorem}\label{thm4.6}
If $G$ is a cycle, then $F_+(G) = 1$.  
\end{theorem}
\proof
Suppose $G = C_n$. Choose a vertex of $G$ and label it as $v_1$. Move clockwise around the cycle and label the vertices as $v_2, \ldots, v_n$. Let $S = \{v_i, v_j\}$ for any distinct $i, j \in \{1, \ldots, n\}$. Color all vertices in $S$ in blue. Without loss of generality, assume $i < j$. If $j=i+1$, then $S$ is composed of adjacent vertices. These vertices are enough to force all the other vertices under the standard color-change rule. If $j > i+1$, then $V(G) \backslash S$ induces a disconnected graph with the paths $P_{j-i-1}$ and $P_{(n-2)-(j-1)}$ as connected components. In each path, $v_i$ and $v_j$ are endpoints. By an argument in the proof of \ref{thm4.5}, the white vertices can now be forced. Hence, $S$ is positive semidefinite zero forcing set for any set of two vertices of the graph. Thus, $F_+(G) \geq 1$. Since a single vertex in $C_n$ is failed, then $F_+(G) = 1$. 

The next proposition and corollary are analogous to Proposition \ref{prop3.2} and Corollary \ref{cor3.3}. The proofs are similar except that a different color-change rule is used.

\begin{proposition}\label{prop4.7}
Suppose $G$ is a disconnected graph. Let $G_1,  \ldots , G_k$ be its connected components. For all $i \in \{1,  \ldots , k\}$, the following holds: If $F_i$ is a maximum failed positive semidefinite zero forcing set for $G_i$, then the set $V(G) \backslash (V(G_i) \backslash F_i)$ is a maximal failed positive semidefinite zero forcing set for $G$.  
\end{proposition}

\begin{corollary}\label{cor4.8}
Suppose $G$ is a disconnected graph. Let $G_1,  \ldots , G_k$ be its connected components. Then $F_+(G) = |G| -$ min$\{|G_i| - F_+(G_i)$ $:$ $i = 1,  \ldots , k\}$. 
\end{corollary}


\begin{proposition}\label{prop4.10} 
Let $G$ be a graph on $n$ vertices and $X$ be a module of order $k$, where $k > 1$. If $X$ contains a connected component with more than one vertex, then $F_+(G) \geq n-k$. 
\end{proposition} 
\proof
It is enough to show that the proposition holds for the case where two vertices of $X$ are adjacent. Let $u$ and $v$ be these vertices. Let $S = V(G) \backslash X$. Color all vertices in $S$ blue. Since $X$ is a module, $u$ and $v$ have the same neighbors in $S$. Forcing under the positive semidefinite color-change rule may only occur if $u$ and $v$ are not adjacent, which is not the case. Hence $S$ is failed and so $F_+(G) \geq |S| = | V(G) \backslash X| =  n - |X| = n-k$.
\qed

\begin{theorem}\label{thm4.12} 
Let $G$ be a connected graph with $n$ vertices. Then $F_+(G)=n-2$ if and only if $G$ has a module composed of two adjacent vertices. 
\end{theorem} 
\proof
Let $S = V(G) \backslash X$, where $X$ is a module composed of two adjacent vertices. Since the vertices in $X$ are adjacent, the induced subgraph of $X$ is a path on two vertices. By Proposition \ref{prop4.10}, $F_+(G) \geq n-2$. Since $G$ is connected, by Proposition \ref{prop4.1} and Theorem \ref{thm4.2}, $F_+(G) < n-1$. Thus, $F_+(G) = n-2$. Conversely, suppose $F_+(G)=n-2$. Let $S$ be a failed positive semidefinite zero forcing set with $|S| = n-2$. If $V(G) \backslash S$ is not connected, then $S$ forces the two vertices, say $u$ and $v$, not in $S$, a contradiction. Hence, $uv \in E(G)$. If there exists a neighbor $w$ of $u$ that is not adjacent to $v$, then $w \rightarrow u$, a contradiction. Therefore $u$ and $v$ must have the same set of neighbors. This implies that $u$ and $v$ form a module of order $2$ of adjacent vertices.
\qed

Consider the complete graph with $n$ vertices, where $n \geq 2$. Then any pair of vertices in $K_n$ forms a module of order $2$ and the vertices in the module are adjacent. Hence, the next corollary follows directly from Theorem \ref{thm4.12}.

\begin{corollary}\label{cor4.13} 
If $n \geq 2$, then $F_+(K_n) = n-2$.
\end{corollary} 

\indent We now look at the case where $F_+(G)$ is very low. That is, $F_+(G) = 0$ or $F_+(G) = 1$. Note that if $F_{+}(G) =0$, then $Z_{+}(G)=1$ and any vertex in $V(G)$ is a positive semidefinite zero forcing set. It is natural to ask if the converse holds. 

\begin{proposition}\label{prop4.14}
If $F_+(G) = 0$, then $G$ is connected. 
\end{proposition}
\proof
Let $F_{+}(G) = 0$. Suppose $G$ has more than one connected components. Choose one of the connected components of $G$, say $W$. Color all the vertices in $S = V(G) \backslash W$ blue. Since $S$ and $W$ are disconnected, $S$ is failed. Since $G$ is disconnected with $W$ as one of its connected components $0 < |W| < |V(G)|$. Thus, $F_+(G) \geq |S| = |V(G) \backslash W| = |V(G)| - |W| > |V(G)| - |V(G)| = 0$, a contradiction. Therefore, $G$ is connected. 
\qed

\begin{proposition}\label{prop4.15}
If $F_{+}(G)=0$, then $G$ has no cycles.
\end{proposition}
\proof
Suppose $G$ has a cycle. Consider the biggest cycle $C_m$ contained in $G$ as a subgraph (possibly not induced). Then, $m$ must be at least 3. Starting from a chosen vertex, move counterclockwise along $C_m$ and label its vertices as $v_1, v_2,  \ldots  v_m$. Consider $S = \{v_2\}$. Color the vertex in $S$ blue. Since $v_2$ is adjacent to $v_1$ and $v_3$, both of which are not in $S$, and $G[V(C_m) \backslash \{v_2\}]$ is connected, we must have that $S$ is stalled. Therefore, $F_+(G) \geq |S| = 1 > 0$.  
\qed

\begin{theorem}\label{thm4.16}
Let $G$ be a graph. Then $F_+(G) = 0$ if and only if $G$ is a tree. 
\end{theorem}
\proof
 Suppose $F_+(G) = 0$. By Propositions \ref{prop4.14} and \ref{prop4.15}, $G$ is connected and has no cycles. Therefore, $G$ is a tree. \\
\indent Conversely, suppose $G$ is a tree. If $v$ is a leaf, then it forces its neighbor. If $v$ is not a leaf, we can decompose $G$ into $W_1$,  \ldots , $W_k$, and $\{v\}$, where the $W_i$'s are the mutually disjoint connected components of $G - v$. Since $G$ is a tree, $G$ has no cycles, which implies that $G[\{v\} \cup W_i]$ is a tree (smaller than $G$) for $i = 1,  \ldots , k$. For each induced tree, $v$ is a leaf and so it can force one of the vertices in each $W_i$. Now, let $w_i \in V(W_i)$ be the forced vertex for each $i$. We apply to $w_i$ what we did to $v$. We repeat these steps and observe that decomposition takes place where smaller trees are obtained. This process must terminate because our graph is finite. The last iteration involves forcing the leaves of $G$. Thus, every singleton in $V(G)$ is a positive semidefinite zero forcing set for $G$. We conclude that $F_+(G)=0$. 
\qed

\begin{corollary}\label{cor4.17}
Let $G$ be a graph. Then $F_+(G) = 0$ if and only if $Z_+(G)=1$. 
\end{corollary}
\proof
Suppose $F_+(G) = 0$. Since $F_+(G)$ is the maximum cardinality of any failed positive semidefinite zero forcing set for $G$, any singleton in $V(G)$ is not failed. Hence, $Z_+(G) = 1$. Conversely, suppose that $Z_+(G) = 1$. By \cite[Theorem $3.5$]{Barioli}, $M_+(G) \leq Z_+(G) = 1$.
By \cite[Theorem $4.1$]{HVDH}, $G$ is a tree. Finally, by Theorem \ref{thm4.16}, $F_+(G) = 0$. 
\qed

\indent We next characterize graphs $G$ for which $F_+(G) = 1$. To aid us in our proof, we first remark that given a connected graph $G$, if $Z_+(G) = 2$, then exactly one block of $G$ must have a cycle \cite{Ekstrand}.

\begin{theorem}\label{thm4.18}
Let $G$ be a graph. Then $F_+(G) = 1$ if and only if $G$ is a set of two isolated vertices or $G$ is a cycle. 
\end{theorem}
\proof
 If $G$ is a set of two isolated vertices, then by Theorem \ref{thm4.2}, $F_+(G) = n - 1 = 2 - 1 = 1$. If $G$ is a cycle, by Theorem \ref{thm4.6}, $F_+(G) = 1$. \\
\indent Conversely, let $F_+(G) = 1$. Suppose that $G$ is a disconnected graph. Let $G_1,  \ldots  , G_k$ be its connected components. Clearly, $F_+(G) \geq |G_i|$ for each $i$. Then, $1 \geq |G_i|$. This implies that $|G_i| = 1$ for all $i \in \{1,  \ldots , k\}$. Therefore, each $G_i$ is an isolated vertex. This implies that $F_+(G_i) = |G_i| - 1 = 0$ for all $i \in \{1,  \ldots , k\}$. By Corollary \ref{cor4.8}, $F_+(G) = |G| - min\{|G_i| - F_+(G_i)\}$. Thus, $1$ $=$ $|G| - (|G_1| + F(G_1)$ $=$ $|G| - (1 + 0)$. Finally, we get $|G| = 2$. Therefore, $G$ must be a set of two isolated vertices. \\
\indent Now, suppose $G$ is connected. Suppose $F_+(G) = 1$. Note that $Z_+(G) - 1 \leq F_+(G) = 1$. Hence, $Z_+(G) \geq 2$. By Theorem \ref{thm4.16}, $G$ is not a tree, and so by Corollary \ref{cor4.17}, $Z_+(G) > 1$. Therefore $Z_+(G) = 2$. From \cite[Theorem $4.1$]{Ekstrand}, exactly one block of $G$ must have a cycle. Hence, $G$ is a cycle, or $G$ is formed by adding edges in a cycle, or $G$ has at least two blocks with exactly one block containing a cycle. We show that in the two latter cases, we contradict the fact that $F_+(G) = 1$. This will then imply that $G$ must be a cycle. \\ 
\indent Case 1: Suppose $G$ is obtained from $C_m$ by adding edges to $C_m$, where $m$ is the maximum possible number of vertices in a cycle contained in $G$. Since $G$ properly contains $C_m$, $m > 3$. Label the vertices of $C_m$ in order (counterclockwise starting from one random vertex) as $v_1$,  \ldots  , $v_m$. Let $v_iv_j \in E(G)$ $\backslash$ $E(C_m)$. Without loss of generality, assume $v_i = v_2$. Otherwise, we can just relabel the vertices so that we start from $v_{i-1}$. Consider the set $\{v_1, v_3\}$ and color its elements blue. Observe that $G[V(G) \backslash \{v_1, v_3\}]$ is connected. Therefore, the first iteration of the positive semidefinite color-change rule is just the same as the standard color-change rule. Note that $v_1$ and $v_3$ are both adjacent to $v_2$. Furthermore, $v_1$ is adjacent to $v_m$ and $v_3$ is adjacent to $v_4$. Hence, the set $\{v_1, v_3\}$ is stalled. Therefore $F_+(G) \geq |\{v_1, v_3\}| = 2 > 1.$ \\ 
\indent Case 2: Suppose $G$ has at least two blocks with exactly one block containing a cycle. This implies that exactly one block of $G$ contains a cycle $C_m$ while the remaining blocks are paths on two vertices. Since $G$ is connected, there exists at least one vertex $u$ in a $P_2$-block that is adjacent to a vertex $v$ in the cycle block. Color $u$ and $v$ blue. Consider a vertex $w$ in $V(C_m) \backslash \{v\}$. To show that $S = \{u, v\}$ is a failed positive semidefinite zero forcing set, it suffices to show that $w$ will never be forced in any iteration of the failed positive semidefinite color-change rule. Suppose, for a contradiction, that $w$ can be forced by a vertex $w_1$. \\ 
\indent If $w_1$ belongs to another block of the graph, then $w_1$ has to be in a $P_2$-block. If $w_1 = u$, then $G[V(C_m) \cup \{u, w\}]$ contains a cycle bigger than $C_m$, a contradiction to $C_m$ being of maximum size. Therefore, $w_1 \neq u$. Hence, $w_1 \notin S$. For it to force another vertex, $w_1$ must first be forced. Hence, there must exist $w_2$ which forces $w_1$. If $w_2 = u$, then $G[V(C_m) \cup \{u, w, w_1\}]$ contains a cycle bigger than $C_m$, a contradiction. Hence, $w_2 \neq u$. We repeat the process. Since our graph is finite, it must terminate. Hence, there does not exist a $w_j$ belonging to another block for which $w_j \rightarrow w_{j-1} \rightarrow  \ldots  \rightarrow w_1 \rightarrow u$. \\ 
\indent Suppose $w_1$ belongs to the block containing $C_m$. Then $w_1 \neq v$ since $v$ is also adjacent to another vertex in $V(C_m) \backslash \{w_1\}$ and the graph $G[V(C_m) \backslash \{v\}]$ is connected. Therefore, $w_1 \notin S$. Hence, for it to force another vertex, $w_1$ must first be forced by a vertex $x_1$. We repeat the process in the preceding paragraph to conclude that $x_1$ must not belong to a $P_2$-block of $G$. Hence, repeating the same process, a vertex $w \notin S$ in the cycle block can only be forced by a vertex in the cycle-block except from $v$. Since $w$ is arbitrary, no other vertex in the cycle block can be forced. Therefore, $S$ is stalled and thus, a failed positive semidefinite zero forcing set. Hence, for the second case, $F_+(G) \geq |S| = 2 > 1$. This is a contradiction. \\ 
\indent Therefore, if $G$ is connected, it has to be a cycle.
\qed

We already know from the preceding section that $F_+(G) = 0$ when $G$ is a tree. In particular, $F_+(P_n) = 0$. We also know from Theorem \ref{thm4.6} and Corollary \ref{cor4.13} that $F_+(C_n) = 1$ and $F_+(K_n) = n - 2$. We now consider other graph families. It is worth noting that in some instances, it will be helpful to compute first the value of $F(G)$ for a graph $G$. This is because of our next theorem. 

\begin{theorem}\label{thm4.19}
For any graph $G$, $F_+(G) \leq F(G)$. 
\end{theorem}
\proof
Let $S$ be a maximum failed positive semidefinite zero forcing set for $G$. Thus, $S$ is stalled under the positive semidefinite-color change rule. It follows that $S$ is also stalled under the standard color-change rule. Therefore, $F(G) \geq |S| = F_+(G)$. 
\qed

\begin{theorem}\label{thm4.20}
Let $n$ be an integer such that $n \geq 4$. If $W_n$ is the wheel graph with $n$ vertices, then $F_+(W_n) = \lfloor \frac{2n-2}{3} \rfloor$. 
\end{theorem}
\proof 
Let $n \geq 4$. We first prove the case where $n=5$. From Theorem \ref{thm3.6}, we have $F(W_5) = 3$. Thus, by Theorem \ref{thm4.19}, $F_+(W_5) \leq 3$. Note that $W_5$ has no module of two adjacent vertices. Also, $W_5$ is a connected graph that is not a cycle but contains a cycle. Hence, by Theorems \ref{thm4.16}, \ref{thm4.18}, and \ref{thm4.12}, $F_+(W_5) > 1$ and $F_+(W_5) \neq 5 - 2 = 3$. Therefore, $F_+(W_5) = 2$. Since $\lfloor \frac{2(5)-2}{3} \rfloor = \lfloor \frac{8}{3} \rfloor = 2$, we have $F_+(W_n) = \lfloor \frac{2n-2}{3} \rfloor$ for $n=5$. \\
\indent Now, suppose $n \neq 5$. From Theorem \ref{thm3.6}, we have $F(W_n) = \lfloor \frac{2n-2}{3} \rfloor$. We show that $F_+(W_n) = F(W_n)$. By Theorem \ref{thm4.19}, $F_+(W_n) \leq F(W_n)$. Thus, it remains to show that there exists a failed positive semidefinite zero forcing set with cardinality $\lfloor \frac{2n-2}{3} \rfloor$. Express $n-1$ as $3k$, $3k+1$, or $3k+2$ based on the reminder when $n-1$ is divided by $3$. Label the vertices in the cycle $C_{n-1}$ consecutively as $v_1,  \ldots , v_{n-1}$ starting from one chosen vertex. Label the hub of the wheel as $x$. Denote $v_n$ as $v_1$ and $v_0$ as $v_{n-1}$. \\
\indent Case 1: $n-1 = 3k$. \\
\indent Color the vertices in the following set blue: $S$ $=$ $\{v_{i}$ $:$ $i < n$ and $i \not\equiv 0$ mod $3\}$. If $i \equiv 1$ mod $3$, then the blue vertex $v_{i}$ is adjacent to $x$ and $v_{i-1}$, which are both white. If $i \equiv 2$ mod $3$, then $v_{i}$ is adjacent to $x$ and $v_{i+1}$, which are both white. Furthermore, $G[V(W_n) \backslash S]$ is connected. Hence, $S$ is a failed positive semidefinite zero forcing set. Therefore, $F_+(W_n) \geq |S|$ $=$ $2(\frac{n-1}{3})$ $=$ $\frac{2n-2}{3}$ $=$ $\lfloor \frac{2n-2}{3} \rfloor$. \\ 
\indent Case 2: $n-1 = 3k+1$. \\
\indent Color the vertices in the following set blue: $S$ $=$ $\{v_{i}$ $:$ $i < n-1$ and $i \not\equiv 0$ mod $3\}$. If $i \equiv 1$ mod $3$, then $v_{i}$ is adjacent to $x$ and $v_{i-1}$, which are both white. If $i \equiv 2$ mod $3$, then $v_{i}$ is adjacent to $x$ and $v_{i+1}$, which are both white. Furthermore, $G[V(W_n) \backslash S]$ is connected. Hence, $S$ is a failed positive semidefinite zero forcing set. Therefore, $F_+(W_n) \geq |S|$ $=$ $2(\frac{n-2}{3})$ $=$ $\frac{2n-4}{3}$ $=$ $\frac{2n-2}{3} - \frac{2}{3}$ $=$ $\lfloor \frac{2n-2}{3} \rfloor$. \\
\indent Case 3: $n-1 = 3k+2$. \\
\indent Color the vertices in the following set blue: $S$ $=$ $\{v_{n-2}\}$ $\cup$ $\{v_{i}$ $:$ $i < n-2$ and $i \not\equiv 0$ mod $3\}$. The blue vertex $v_{n-2}$ is adjacent to $x$, $v_{n-1}$ and $v_{n-3}$, all of which are white. Suppose $i \neq n-2$. If $i \equiv 1$ mod $3$, then $v_{i}$ is adjacent to $x$ and $v_{i-1}$, which are both white. If $i \equiv 2$ mod $3$, then $v_{i}$ is adjacent to $x$ and $v_{i+1}$, which are both white. Furthermore, $G[V(W_n) \backslash S]$ is connected. Hence, $S$ is a failed positive semidefinite zero forcing set. Therefore, $F_+(W_n) \geq |S|$ $=$ $2(\frac{n-3}{3}) + 1$ $=$ $\frac{2n-3}{3}$ $=$ $\frac{2n-2}{3} - \frac{1}{3}$ $=$ $\lfloor \frac{2n-2}{3} \rfloor$. \\
\indent In all cases, $F_+(W_n) = \lfloor \frac{2n-2}{3} \rfloor$.
\qed

\begin{theorem}\label{thm4.21}
Let $m$, $n$ $\in$ $\mathbb{N}$. If $p$ $=$ min$\{m, n\}$, then
$$F_+(K_{m,n}) = \begin{cases}
0, & \; if \; p = 1 \\
m + n - 3, & \; if \; p = 2 \\
m + n - 4, & \; if \; p > 2 \\
\end{cases}$$
\end{theorem}

\proof
Let $m$, $n$ $\in$ $\mathbb{N}$ and $p$ $=$ min$\{m, n\}$. If $p = 1$, $K_{m,n}$ is a tree. By Theorem \ref{thm4.16}, $F_+(K_{m,n})$ $=$ $0$. Let $p = 2$. Now, let $S_1$ be the first partite component of $K_{m,n}$ which contains $m$ vertices and $S_2$ be the second partite component of $K_{m,n}$ which contains $n$ vertices. Note that $G[S_1]$ and $G[S_2]$ is a graph of $m$ and $n$ isolated vertices, respectively. Label the vertices of $S_1$ as $a_1,  \ldots , a_m$ and the vertices of $S_2$ as $b_1,  \ldots , b_m$. Consider the set $S = V(K_{m,n}) \backslash S'$, where $S' = \{a_1, b_1, b_2\}$ if $n=2$, and $S' = V(K_{m,n}) \backslash \{a_1, a_2, b_1\}$ if $n \neq 2$. Observe that $G[S']$ is a path on $3$ vertices and thus connected. Therefore, an iteration of the positive semidefinite color-change rule applied to $S$ is simply a standard color-change rule. Observe that the vertices contained in $S$ are all vertices from the same partite component. Since $p=2$, each vertex in $S$ has two white neighbors coming from the other partite component of $K_{m,n}$. Hence, $S$ is failed. Therefore, $F_+(K_{m,n}) \geq |S| = m + n - 3$ when $p = 2$. Observe that $K_{m,n}$ has a module of order $2$. However, the vertices in each module are non-adjacent. Hence, $K_{m,n}$ contains no module of order $2$ of adjacent vertices. Since $K_{m,n}$ is also connected, by Proposition \ref{prop4.1}, Theorem \ref{thm4.2}, and Theorem \ref{thm4.12}, $F_+(K_{m,n}) < m + n - 2$. Therefore, for $p=2$, $F_+(K_{m,n}) = m+n-3$. \\
\indent Lastly, let $p > 2$. Let $S$ contain all vertices of $K_{m,n}$ except two from the first partite and two from the second partite. Let $u_1$ and $u_2$ be the vertices from the first partite that are not contained in $S$, and let $v_1$ and $v_2$ be the vertices from the second partite that are not contained in $S$. Observe that $G[\{u_1, u_2, v_1, v_2\}]$ is a cycle on $4$ vertices and thus connected. Therefore, an iteration of the positive semidefinite color-change rule applied to $S$ is simply a standard color-change rule. Notice that each vertex from the first partite contained in $S$ has $v_1$ and $v_2$ as its white neighbors while each vertex from the second partite contained in $S$ has $u_1$ and $u_2$ as its white neighbors. Therefore, $S$ is failed, implying that $F_+(K_{m,n}) \geq |S| = n+m-4$ for $p > 2$. \\
\indent We next show that for $p > 2$, a set $\overline{S}$ containing $n+m-3$ vertices is a positive semidefinite zero forcing set for $K_{m,n}$. Let $w_1, w_2, w_3$ be the vertices of $K_{m,n}$ not contained in $\overline{S}$. \\
\indent If $w_1, w_2$, and $w_3$ are from the first partite, then $G[V(G) \backslash \overline{S}]$ is a graph on isolated vertices. Hence, each vertex outside $\overline{S}$ can be forced, implying that $\overline{S}$ is a positive semidefinite zero forcing set. \\
\indent Now, without loss of generality, let $w_1$ and $w_2$ be from the first partite and $w_3$ be from the second partite. Since $p > 2$, $w_3$ can be forced by any vertex from the first partite other than $w_1$ and $w_2$. Suppose, for a contradiction, that at least one of $w_1$ and $w_2$ will not be forced no matter how many times we apply the positive semidefinite color-change rule. Then, $\overline{S}$ $\cup$ $\{w_3\}$ must be failed. This implies that $F_+(K_{m,n})$ $\geq$ $|\overline{S}$ $\cup$ $\{w_3\}|$ $=$ $n+m-2$. This is a contradiction since $K_{m,n}$ is connected and has no module of order 2 of adjacent vertices. Hence, $w_1$ and $w_2$ must be eventually forced, implying that $\overline{S}$ is a positive semidefinite zero forcing set. \\
\indent Therefore, for $p > 2$, $F_+(K_{m,n}) = n + m - 4$. 
\qed

\begin{theorem}\label{thm4.22}
Let $G$ be the hypercube $Q_n$. Then $F_+(Q_n) \geq 2^n - n -1$.
\end{theorem}
\proof
Note that $Q_1 = P_2$ and $Q_2 = C_4$. By Theorems \ref{thm4.5} and \ref{thm4.6}, $F_+(Q_1) = 0$ and $F_+(Q_2) = 1$. Suppose $n \geq 3$. Consider the set $F = V(G) \backslash W$ as in the proof of Theorem \ref{thm3.7}. Since $F$ is a failed zero forcing set for $Q_n$, by observation \ref{obs2.4}, any subset of $F$ is also a failed zero forcing set. In particular, $F' = F \backslash \{(1, \ldots, 1)\}$ is a failed zero forcing set. Color all vertices in $F'$ blue and subject $F'$ under the positive semidefinite color-change rule. The white vertices are those vertices in $W \cup \{(1, \ldots, 1)\}$. Since each vertex in $W$ is of the form $(x_1, \ldots, x_n)$ with exactly $n-1$ of the $x_i$'s equal to 1, then each vertex in $W$ is adjacent to $(1, \ldots, 1)$. Hence, $W \cup \{(1, \ldots, 1)\}$ induces a connected subgraph and the first iteration of the positive semidefinite color-change rule is just a standard one. Since $F' = F \backslash \{(1, \ldots, 1)\} = V(G) \backslash (W \cup \{(1, \ldots, 1)\})$ is failed under the standard color-change rule, then $F'$ is stalled. Hence, $F'$ is also a failed positive semidefinite zero forcing set and for $n \geq 3$, $F_+(Q_n) \geq |F| = |Q_n| - |W| - 1 = 2^n - {{n}\choose{n-1}} - 1 = 2^n - n -1 $.
\qed

\begin{theorem}\label{thm4.23}
Let $G$ be the half-graph $H_s$. Then $F_+(H_1) = 0$ and $F_+(H_s) = 2s - 4$ for $s \geq 2$. 
\end{theorem}
\proof
Note that $H_1 = P_2$. By Theorem \ref{thm4.5}, $F_+(H_1) = 0$. Suppose $s \geq 2$. Partition $H_s$ into $S_1$ and $S_2$ such that this partition defines the bipartite property of $H_s$. Label the vertices in $S_1$ as $v_1$,  \ldots , $v_s$ and the corresponding vertices in $S_2$ as $v_{s+1},  \ldots , v_{2s}$, respectively. Let $V = \{v_{s-1}, v_{s}, v_{2s-1}, v_{2s}\}$. The subgraph of $H_s$ induced by $V$ is a path and thus, connected. Thus, the first iteration of the positive semidefinite color-change rule applied to the set $V(H_s)$ $\setminus$ $V$ is just a standard one. However, $v_{s-1}$ cannot be forced since its neighbors $v_{2s-1}$ and $v_{2s}$ are both outside $V(H_s)$ $\setminus$ $V$. Also, $v_{s}$ cannot be forced since its only neighbor $v_{2s}$ is outside $V(H_s)$ $\setminus$ $V$. Lastly, every vertex in $V(H_s)$ $\setminus$ $V$ adjacent to $v_{2s-1}$ is also adjacent to $v_{2s}$. Therefore, $v_{2s-1}$ and $v_{2s}$ cannot be forced. This implies that $V(H_s)$ $\setminus$ $V$ is stalled. Therefore, $V(H_s)$ $\setminus$ $V$ is a failed positive semidefinite zero forcing set and $F_+(H_s) \geq 2s - 4$. \\
\indent Let $G = H_s$. We next show that any set $S$ of $2s-3$ vertices is a positive semidefinite zero forcing set for $G$. Let $V_1$ and $V_2$ be the vertices outside $S$ contained in $S_1$ and $S_2$, respectively. Without loss of generality, assume $|V_1| \leq |V_2|$. We have the following cases: (a) $|V_1| = 0$ and $|V_2| = 3$; (b) $|V_1| = 1$ and $|V_2| = 2$. Suppose $|V_1| = 0$ and $|V_2| = 3$. \\
\indent If $|V_1|=0$ and $|V_2| = 3$, then $G[V_2]$ is a subgraph with $3$ isolated vertices, say $u_1$, $u_2$, and $u_3$. The subgraphs $G[S \cup \{u_i\}]$ for $i \in \{1, 2, 3\}$ are connected graphs with exactly one white vertex, which is $u_i$. Hence, under the positive semidefinite color-change rule, vertices $u_1$, $u_2$, and $u_3$ will be forced. Therefore, Hence, $S$ is a positive semidefinite zero forcing set for $H_s$. \\
\indent Next, suppose $|V_1| = 1$ and $|V_2| = 2$. Let $V_1 = \{v_{j}\}$ and $V_2 = \{v_{s+k}, v_{s+l}\}$ for some $j, k, l \in \{1,  \ldots , s\}$. The subgraph $G[V(H_s) \setminus S]$ is either connected or disconnected with at most three connected components. \\
\indent If $G[V(H_s) \setminus S]$ is disconnected with exactly three connected components, then $G[V(H_s) \setminus S]$ is a graph with three isolated vertices. By a previous argument, these vertices can be forced. \\
\indent Suppose $G[V(H_s) \setminus S]$ is disconnected with exactly two connected components. Without loss of generality and since $G[V_2]$ is a graph on isolated vertices, the connected components of $G[V(H_s) \setminus S]$ are $G[\{v_j, v_{s+k}\}]$ and $G[\{v_{s+l}\}]$. It follows from the definition of a half-graph that $l < j \leq k$. Since $G[\{v_{s+l}\}]$ is an isolated vertex, $v_{s+l}$ will be forced. Hence, in an updated coloring of $G$, there exists only two white vertices. Since $G$ is connected and has no module of order 2, the two white vertices will also be forced. \\
\indent Lastly, suppose $G[V(H_s) \setminus S]$ is connected. Then $j \leq k$ and $j \leq l$. Without loss of generality, assume that $j \leq k < l$. Then $v_l \rightarrow v_{s+l} \rightarrow v_{j} \rightarrow v_{s+k}$. \\
\indent In all cases, a set $S$ of $2s-3$ vertices is a positive semidefinite zero forcing set for $G$. Hence, $F_+(H_s) \leq 2s-4$. Therefore, the equality $F_+(H_s) = 2s-4$ holds for $s \geq 2$. 
\qed

\section{Failed zero forcing parameter and the minimum rank problem}\label{sec5}

\indent We first present known results on the maximum nullity as well as the zero forcing number of different graph families.  

\begin{center}
\begin{table}[H]
\begin{tabular}{lllllll}
\hline
Reference 	& $G$ & Order & $M(G)$ & $Z(G)$ & $M_+(G)$ & $Z_+(G)$ \\ 
\hline
\cite{AIM}, \cite{HVDH} 		& $P_n$ & $n$ & 1 & 1 & 1 & 1 \\ 
\cite{AIM}, \cite{HVDH} 		& $C_n$ & $n$ & 2 & 2 & 2 & 2 \\
\cite{AIM}, \cite{FH} 			& $K_n$ & $n$ & $n-1$ & $n-1$ & $n-1$ & $n-1$ \\
\cite{AIM}, \cite{Peters} 	& $Q_1$ & $2$ & 1 & 1 & 1 & 1  \\
\cite{AIM}, \cite{Peters} 	& $Q_2$ & $4$ & 2 & 2 & 2 & 2  \\
\cite{AIM}, \cite{Peters} 	& $Q_n, n \geq 3$ & $2^n$ & $2^{n-1}$ & $2^{n-1}$ & $2^{n-1}$ & $2^{n-1}$ \\
\cite{Peters} 							& $W_n, n \geq 4$& $n$ & 3 & 3 & 3 & 3 \\
\cite{Peters} 							& $K_{m,n}, m \geq n$ & $m+n$ & $m+n-2$ & $m+n-2$ & $n$ & $n$ \\
\cite{AIM}, \cite{Peters} 				& $H_s$ & $2s$ & $s$ &  $s$ & $s$ & $s$ \\
\end{tabular}
\caption{\label{table5.1} Maximum nullities and zero forcing numbers of several graph families}
\end{table}
\end{center}

In \cite{FJS}, graphs for which $F(G)$ is strictly less than $Z(G)$ are characterized. However, for $F_+(G)$, we will only consider the graph families listed in Table \ref{table5.1}

\begin{theorem}\label{thm5.1}
\cite[Theorem $2.8$]{FJS} Let $G$ be a graph. $F(G) < Z(G)$ if and only if $G = K_n$ or $G = \overline{K_n}$.
\end{theorem}

The next theorem is an immediate consequence of comparing the values of $Z_+(G)$ in Table \ref{table5.1} and the results for $F_+(G)$ in Section \ref{sec4}. 

\begin{theorem}\label{thm5.2}
Let $G$ be a graph. 
\begin{enumerate}
\item If $G$ is a path, cycle, complete graph, tree, or a set of isolated vertices, then $F_+(G) < Z_+(G)$. 
\item If $G = Q_n$, then $F_+(G) < Z_+(G)$ if and only if $n = 1$ or $n =2$.
\item If $G = W_n$ with $n \geq 4$, then $F_+(G) < Z_+(G)$ if and only if $n = 4$ or $n = 5$.
\item Let $G = K_{m,n}$. Then $F_+(G) < Z_+(G)$ if and only if min$\{m, n\}$ $= 1$, or $G = K_{2,2}$, or $G=K_{3,3}$.
\item If $G = H_s$, then $F_+(G) < Z_+(G)$ if and only if $s \leq 3$. 
\end{enumerate}
\end{theorem}

\indent We now relate $F(G)$ and $F_+(G)$ with the minimum rank problem. It is known that the zero forcing number is an upper bound for the maximum nullity of a graph. In particular, we have the following theorems:

\begin{theorem}\label{thm5.3}
\cite[Proposition $2.4$]{AIM} For any graph $G$, $M(G) \leq Z(G)$.
\end{theorem}

\begin{theorem}\label{thm5.4}
\cite[Theorem $3.5$]{Barioli} For any graph $G$, $M_+(G) \leq Z_+(G)$.
\end{theorem}

\indent The bounds in Theorems \ref{thm5.3} and \ref{thm5.4} are sharp. The graph families listed in Table \ref{table5.1} all satisfy the equalities $M(G) = Z(G)$ and $M_+(G) = Z_+(G)$. By the rank-nullity theorem, we know that $mr(G) + M(G) = |G|$. Since there are graphs for which $M(G) = Z(G)$ and failed zero forcing sets are sets that are not zero forcing sets, it is natural to ask if there are graphs for which $mr(G) = F(G)$. \\
\indent In spite of an abundant research on the area of zero forcing, a full classification of graphs for which $M(G) = Z(G)$ or $M_+(G) = Z_+(G)$ remains to be an open problem. Hence, we also expect that a full classification of graphs for which $\text{mr}(G) = F(G)$ or $\text{mr}_+(G) = F_+(G)$ is also a difficult problem. In this paper, we will only consider the graph families in Table \ref{table5.1}. \\
\indent The next two propositions directly follow from the rank-nullity theorem.

\begin{proposition}\label{prop5.5}
Suppose $Z(G) = M(G)$. Then $F(G) + Z(G) = |G|$ if and only if $F(G) = mr(G)$.
\end{proposition}

\begin{proposition}\label{prop5.6}
Suppose $Z_+(G) = M_+(G)$. Then $F_+(G) + Z_+(G) = |G|$ if and only if $F_+(G) = mr_+(G)$.
\end{proposition}

\indent For graphs $G$ listed in Table \ref{table5.1}, $Z(G) = M(G)$ and $Z_+(G) = M_+(G)$. Hence, to characterize the graphs under each considered graph family for which $F(G) = mr(G)$ or $F_+(G) = mr_+(G)$, it suffices to characterize the graphs under each graph family for which $F(G) + Z(G) = |G|$ or $F_+(G) + Z_+(G) = |G|$. We accomplish this task by using Table \ref{table5.1} as well as the results from Sections \ref{sec3} and \ref{sec4}. From manual computation, we get the next two theorems. 

\begin{theorem}\label{thm5.7}
Let $G$ be a graph. 
\begin{enumerate}
\item If $G$ is a path, then $F(G) < \text{mr}(G)$.
\item If $G$ is a cycle, then $F(G) = \text{mr}(G)$ if and only if $G = C_3$ or $G = C_4$. 
\item If $G$ is a complete graph, then $F(G) = \text{mr}(G)$ if and only if $G = K_3$. 
\item If $G$ is a hypercube, then $F(G) = \text{mr}(G)$ if and only if $G=Q_2$. 
\item If $G$ is a wheel, then $F(G) = \text{mr}(G)$ if and only if $G = W_6$ or $G=W_7$. 
\item If $G = K_{m,n}$ with $m \geq n$, then $F(G) = \text{mr}(G)$ if and only if $m+n=4$. 
\item If $G$ is a half-graph, then $F(G) = \text{mr}(G)$ if and only if $G = H_s$.
\end{enumerate}
\end{theorem}

\begin{theorem}\label{thm5.8}
Let $G$ be a graph. 
\begin{enumerate}
\item If $G$ is a path, then $F_+(G) < \text{mr}_+(G)$.
\item If $G$ is a cycle, then $F_+(G) = \text{mr}_+(G)$ if and only if $G = C_3$. 
\item If $G$ is a complete graph, then $F_+(G) = \text{mr}_+(G)$ if and only if $G = K_3$. 
\item If $G$ is a hypercube, then $F_+(G) = \text{mr}_+(G)$ if and only if $G=Q_3$. 
\item If $G$ is a wheel, then $F_+(G) = \text{mr}_+(G)$ if and only if $G = W_5$, or $G=W_6$, or $G=W_7$. 
\item If $G = K_{m,n}$ with $m \geq n$, then $F_+(G) = \text{mr}_+(G)$ if and only if $n=4$. 
\item If $G$ is a half-graph, then $F(G) = \text{mr}(G)$ if and only if $G = H_4$.
\end{enumerate}
\end{theorem}

\section*{Acknowledgement}

The work of P. A. B. Pelayo was funded by the Accelerated Science and Technology Human
Resource Development Program - National Science Consortium (ASTHRDP-NSC) of the Department of
Science and Technology - Science Education Institute (DOST-SEI), Philippines.

\end{document}